\documentclass[a4paper,oneside,11pt]{article}
\usepackage{amssymb}
\usepackage{amsfonts}
\usepackage{amsmath}

\input{tcilatex}

\begin{document}

\title{The fixed set of the inverse involution on a Lie group}
\author{Haibao Duan\thanks{%
The first author is supported by 973 Program 2011CB302400 and NSFC 11131008.}
and Shali Liu \and Institute of Mathematics, Chinese Academy of Sciences
\and dhb@math.ac.cn}
\maketitle

\begin{abstract}
In \cite{[DL]} we have determined the isomorphism type of the centralizer of
an element in a simpe Lie group. As a sequel to \cite{[DL]} we present a
general procedure to calculate the isomorphism type of the fixed set $%
Fix(\gamma )$ of the inverse involution $\gamma $ on a Lie group $G$.

\begin{description}
\item \textsl{2000 Mathematical Subject Classification: }55M20,\textsl{\ }%
22E15; 53C35.

\item \textsl{Key words and phrases:} Lie group; fixed point; adjoint orbit
\end{description}
\end{abstract}

\section{Introduction}

Let $G$ be a compact, connected and simple Lie group with group unit $e\in G$%
. \textsl{The inverse involution} on $G$ is the periodic $2$ transformation $%
\gamma $ sending each group element $g\in G$ to its inverse $g^{-1}\in G$.
In this paper we present a general procedure to calculate the isomorphism
type of the fixed set

\begin{quote}
$Fix(\gamma )=\{g\in G\mid g=g^{-1}\}$
\end{quote}

\noindent of the involution $\gamma $.

Given a group element $x\in G$ let $M_{x}$, $C_{x}\subset G$ be the adjoint
orbit through $x$ and the centralizer of $x$ in $G$, respectively. That is

\begin{quote}
$M_{x}=\{gxg^{-1}\in G\mid g\in G\}$; $C_{x}=\{g\in G\mid gx=xg\}$.
\end{quote}

\noindent The map $G\rightarrow G$ by $g\rightarrow gxg^{-1}$ is constant
along the left cosets of $C_{x}$ in $G$, and induces a diffeomorphism from
the homogeneous space $G/C_{x}$ onto the orbit space $M_{x}$

\begin{enumerate}
\item[(1.1)] $f_{x}:G/C_{x}\overset{\cong }{\rightarrow }M_{x}$, $%
[g]\rightarrow gxg^{-1}$.
\end{enumerate}

\noindent In view of this identification the isomorphism type of the orbit
space $M_{x}$ is completely determined by the centralizer $C_{x}$. It is
crucial to notice that $x\in Fix(\gamma )$ implies that $M_{x}\subset
Fix(\gamma )$. Naturally, one asks for a partition of the space $Fix(\gamma
) $ by certain adjoint orbits $M_{x}$, and determine the isomorphism types
of the corresponding centralizers $C_{x}$.

Concerning the applications of our approach we assume the reader's
familiarity with the classification on Lie groups. In particular, all $1$%
--connected compact simple Lie groups consists of the three infinite
families $SU(n+1),Sp(n),$ $Spin(n+2)$, $n\geq 2,$ of \textsl{classical groups%
}, and the five \textsl{exceptional Lie groups} $%
G_{2},F_{4},E_{6},E_{7},E_{8}$. For a classical Lie group $G$ the fixed set $%
Fix(\gamma )$ can be easily calculated using linear algebra, see Frankel 
\cite{[F]}. For this reason we shall restrict ourself to the simple
exceptional Lie groups. Explicitly we shall have

\begin{quote}
$G=G_{2},F_{4},E_{6},E_{7},E_{8}$ or $E_{6}^{\ast },E_{7}^{\ast }$,
\end{quote}

\noindent where $G^{\ast }=G/\mathcal{Z}(G)$ with $\mathcal{Z}(G)$ the
center of $G$.

Fix a maximal torus $T$ in $G$ and let $\exp :L(T)\rightarrow T$ be the
exponential map, where $L(T)$ is the tangent space to $T$ at the unit $e$.
In term of a set $\Omega =\{\omega _{1},\cdots ,\omega _{n}\}\subset L(T)$
of fundamental dominant weights of $G$ (see Definition 2.3), together with
the fundamental Weyl cell $\Delta $ corresponding to $\Omega $, our main
result is stated below.

Let $SO(n)$ and $Ss(n)$ be the \textsl{special orthogonal group} and the 
\textsl{semispinor group} of order $n$, respectively. For a connected Lie
group $H$ write $[H]^{2}$ for the group with two components whose identity
component is $H$. For two manifolds $M$ and $N$ denote by $M\coprod N$ their
disjoint union.

\bigskip

\noindent \textbf{Theorem 1.1.} \textsl{For a simple Lie group }$G$\textsl{\
there is a subset }$\mathcal{F}_{G}\subset $\textsl{\ }$\Delta $\textsl{\ so
that}

\begin{enumerate}
\item[(1.2)] $Fix(\gamma )=\{e\}\coprod\limits_{u\in \mathcal{F}_{G}}M_{\exp
(u)}$\textsl{.}
\end{enumerate}

\noindent \textsl{Moreover, for each exceptional Lie group }$G$ \textsl{the
set }$\mathcal{F}_{G}$\textsl{, as well as the isomorphism type of the
adjoint orbit} $M_{\exp (u)}$\textsl{\ with }$u\in \mathcal{F}_{G}$\textsl{%
,\ is tabulated below}

\begin{center}
\begin{tabular}{l|l|l}
\hline\hline
$G$ & $\mathcal{F}_{G}$ & ${\small M}_{\exp (u)}{\small =G/C}_{\exp (u)}$, $%
u\in \mathcal{F}_{G}$ \\ \hline
$G_{2}$ & $\{\frac{{\small \omega }_{{\small 1}}}{2}\}$ & ${\small G}_{2}%
{\small /SO(4)}$ \\ 
$F_{4}$ & $\{\frac{{\small \omega }_{{\small k}}}{2}\}_{{\small k=1,4}}$ & $%
{\small F}_{4}{\small /Spin(9),}$ ${\small F}_{4}{\small /}\frac{{\small %
Sp(3)\times Sp(1)}}{{\small Z}_{2}}$ \\ 
$E_{6}$ & $\{\frac{{\small \omega }_{{\small 2}}}{2},\frac{\omega
_{1}+\omega _{6}}{2}\}$ & ${\small E}_{6}{\small /}\frac{{\small SU(2)\times
SU(6)}}{{\small Z}_{2}}{\small ,}$ ${\small E}_{6}{\small /}\frac{{\small %
Spin(10)\times S}^{{\small 1}}}{{\small Z}_{4}}$ \\ 
$E_{6}^{\ast }$ & ${\small \{}\frac{{\small \omega }_{{\small k}}}{2}{\small %
\}}_{{\small k=1,2}}$ & ${\small E}_{6}^{\ast }{\small /}\frac{{\small %
Spin(10)\times S}^{{\small 1}}}{{\small Z}_{4}}{\small ,}$ ${\small E}%
_{6}^{\ast }{\small /}\frac{{\small SU(2)\times }\frac{{\small SU(6)}}{Z_{3}}%
}{{\small Z}_{2}}$ \\ 
$E_{7}$ & $\{\frac{{\small \omega }_{{\small k}}}{2},\omega _{{\small 7}}\}_{%
{\small k=1,6}}$ & ${\small E}_{7}{\small /}\frac{{\small Spin(12)\times
SU(2)}}{{\small Z}_{2}}{\small ,}$ ${\small E}_{7}{\small /}\frac{{\small %
Spin(12)\times SU(2)}}{{\small Z}_{2}}{\small ,}${\small \ }${\small \exp
(\omega _{7})}$ \\ 
$E_{7}^{\ast }$ & ${\small \{}\frac{{\small \omega }_{{\small k}}}{2}{\small %
\}}_{{\small k=1,2,7}}$ & ${\small E}_{7}^{\ast }{\small /}\frac{{\small %
Ss(12)\times SU(2)}}{{\small Z}_{2}}{\small ,}$ ${\small E}_{7}^{\ast }%
{\small /}[\frac{{\small SU(8)}}{Z_{4}}{\small ]}^{2}{\small ,}$ ${\small E}%
_{7}^{\ast }{\small /}[\frac{{\small E}_{6}{\small \times S}^{1}}{{\small Z}%
_{3}}]^{2}$ \\ 
$E_{8}$ & $\{\frac{{\small \omega }_{{\small k}}}{2}\}_{{\small k=1,8}}$ & $%
{\small E}_{8}{\small /Ss(16),}$ ${\small E}_{8}{\small /}\frac{{\small E}_{%
{\small 7}}{\small \times SU(2)}}{{\small Z}_{2}}$ \\ \hline\hline
\end{tabular}

{\small Table 1. The fixed sets of the inverse involution on exceptional Lie
groups}
\end{center}

Historically, the problem of determining the isomorphism type of the fixed
set $Fix(\gamma )$ of a simple Lie group $G$ has been studied by Frankel 
\cite{[F]} for the classical Lie groups, and by Chen, Nagano \cite{[CN1],[N]}%
, Yokota \cite{[Y1],[Y2],[Y3]} for the exceptional Lie groups, see Remark
4.4. These works rely largely on the specialities of each individual Lie
group and the calculations were performed case by case. In comparison, our
approach is free of the types of simple Lie groups, and is ready to extend
to general cases, see Corollaries 5.1--5.2 of Section 5.

The paper is arranged as follows. Section \S 2 contains a brief introduction
to the roots and weight systems of simple Lie groups. In Section \S 3 the
set $\mathcal{F}_{G}$ specifying the partition on $Fix(\gamma )$ in formula
(1.2) is largely determined by Lemma 3.3. Combining Lemma 3.3 with the
algorithm calculating the isomorphism type of a centralizer $C_{x}$ obtained
in \cite{[DL]}, Theorem 1.1 is established in Section 4. Finally, general
structure of the fixed set $Fix(\gamma )$ of the inverse involution $\gamma $
on an arbitrary Lie group $G$ is discussed briefly in Section 5.

\section{Geometry of roots and weights}

For a simple Lie group $G$ with Lie algebra $L(G)$ and a maximal torus $T$,
the dimension $n=\dim T$ is called the \textsl{rank} of $G$, and the
subspace $L(T)$ of $L(G)$ is called the\textsl{\ Cartan subalgebra} of $G$.
Equip $L(G)$ with an inner product $(,)$ so that the adjoint representation
acts on $L(G)$ as isometries, and let

\begin{quote}
$d:G\times G\rightarrow \mathbb{R}$ (resp. $d:T\times T\rightarrow \mathbb{R}
$)
\end{quote}

\noindent be the induced metric on $G$ (resp. on $T$).

The restriction of the exponential map $\exp :L(G)\rightarrow G$ to $L(T)$
defines a set $\mathcal{S}(G)=\{L_{1},\cdots ,L_{m}\}$ of $m=\frac{1}{2}%
(\dim G-n)$ hyperplanes in $L(T)$, namely, the set of\textsl{\ singular
hyperplanes }through the origin in $L(T)$ \cite[p.168]{[BD]}. Let $%
l_{k}\subset L(T)$ be the normal line of the plane $L_{k}$ through the
origin, $1\leq k\leq m$. Then the map $\exp $ carries $l_{k}$ onto a circle
subgroup of $G$.

\bigskip

\noindent \textbf{Definition 2.1.} Let $\pm \alpha _{k}\in l_{k}$ be the
non--zero vectors with minimal length so that $\exp (\pm \alpha _{k})=e$, $%
1\leq k\leq m$. The subset

\begin{quote}
$\Phi =\{\pm \alpha _{k}\in L(T)\mid 1\leq k\leq m\}$
\end{quote}

\noindent of $L(T)$ is called the \textsl{root system of }$G$.

The \textsl{Weyl group} of $G$, denoted by $W$, is the subgroup of $%
Aut(L(T)) $ generated by the reflections $r_{k}$ in the hyperplane $L_{k}$, $%
1\leq k\leq m$. $\square $

\bigskip

\noindent \textbf{Remark 2.2. }We point out that\textbf{\ }the root system $%
\Phi $ by Definition 2.1 is \textsl{dual} to those that are commonly used in
literatures, e.g. \cite{[B],[Hu]}. In particular, the symplectic group $%
Sp(n) $ is of the type $B_{n}$, while the spinor group $Spin(2n+1)$ is of
the type $C_{n}$.$\square $

\bigskip

The planes in $\mathcal{S}(G)$ divide $L(T)$ into finitely many convex open
cones, called the \textsl{Weyl chambers} of $G$. Fix once and for all a
regular point $x_{0}\in L(T)\backslash \underset{1\leq k\leq m}{\cup }L_{m}$%
, and let $\mathcal{F}(x_{0})$ be the closure of the Weyl chamber containing 
$x_{0}$. Assume that $L(x_{0})=\{L_{1},\cdots ,L_{n}\}$ is the subset of $%
\mathcal{S}(G)$ consisting of the walls of $\mathcal{F}(x_{0})$, and let $%
\alpha _{i}\in \Phi $ be the root normal to the wall $L_{i}\in $ $L(x_{0})$
and pointing toward $x_{0}$. Then the subset $\Delta (x_{0})=\{\alpha
_{1},\cdots ,\alpha _{n}\}$ of $\Phi $ is called \textsl{the} \textsl{system
of} \textsl{simple roots }of $G$ relative to $x_{0}$.

\bigskip

\noindent \textbf{Definition 2.3 }(\cite[p.67]{[Hu]})\textbf{.} Each root $%
\alpha \in \Phi $ gives rise to a linear map

\begin{enumerate}
\item[(2.1)] $\alpha ^{\ast }:L(T)\rightarrow \mathbb{R}$ by $\alpha ^{\ast
}(x)=2(x,\alpha )/(\alpha ,\alpha )$,
\end{enumerate}

\noindent called \textsl{the inverse root} of $\alpha $. The \textsl{weight
lattice of} $G$ is the subset of $L(T)$

\begin{quote}
$\Lambda =\{x\in L(T)\mid \alpha ^{\ast }(x)\in \mathbb{Z}$ for all $\alpha
\in \Phi \}$,
\end{quote}

\noindent whose elements are called \textsl{weights}. Elements in the subset
of $\Lambda $

\begin{enumerate}
\item[(2.2)] $\Omega =\{\omega _{i}\in L(T)\mid \alpha _{j}^{\ast }(\omega
_{i})=\delta _{i,j},$ $\alpha _{j}\in \Delta (x_{0})\}$
\end{enumerate}

\noindent are called the \textsl{fundamental dominant weights} of $G$
relative to $x_{0}$, where $\delta _{i,j}$ is the Kronecker symbol.$\square $

\bigskip

To be precise we adopt the convention that for each simple group $G$ with
rank $n$ its fundamental dominant weights $\omega _{1},\cdots ,\omega _{n}$
are ordered by the order of their corresponding simple roots pictured as the
vertices in the Dynkin diagram of $G$ in \cite[p.58]{[Hu]}. Useful
properties of the weights are:

\bigskip

\noindent \textbf{Lemma 2.4. }\textsl{Let\ }$\Omega =\{\omega _{1},\cdots
,\omega _{n}\}$\textsl{\ be the set of fundamental dominant weights relative
to the regular point }$x_{0}$\textsl{. Then}

\textsl{i) }$\Omega =\{\omega _{1},\cdots ,\omega _{n}\}$\textsl{\ is a
basis for }$\Lambda $ \textsl{over} $\mathbb{Z}$\textsl{;}

\textsl{ii) for each }$1\leq i\leq n$\textsl{\ the half line }$\{t\omega
_{i}\in L(T)\mid t\in \mathbb{R}^{+}\}$\textsl{\ is the edge of the Weyl
chamber }$\mathcal{F}(x_{0})$\textsl{\ opposite to the wall }$L_{i}$\textsl{;%
}

\textsl{iii) if }$G$\textsl{\ is simple, }$(\omega _{i},\omega _{j})>0$ 
\textsl{for all }$1\leq i,j\leq n$\textsl{.}$\square $

\noindent \textbf{Proof.} Property i) is well known. By (2.2) each weight $%
\omega _{i}\in \Omega $ is perpendicular to $\alpha _{j}$ (i.e. $\omega
_{i}\in L_{j}$) for all $j\neq i$, $1\leq j\leq n$. This verifies ii). For
iii) we refer to \cite[p.72, Exercise 8]{[Hu]}.$\square $

\bigskip

Let $\mathcal{Z}(G)$ be the center of the group $G$, and let $\Lambda
_{e}=\exp ^{-1}(e)\subset L(T)$ be the \textsl{unit lattice}. The set $%
\Delta (x_{0})=\{\alpha _{1},\cdots ,\alpha _{n}\}$ of simple roots spans
also a lattice $\Lambda _{r}$ on $L(T)$, known as the \textsl{root lattice }%
of $G$.

\bigskip

\noindent \textbf{Lemma 2.5 (\cite[(3.3)]{[DL]}) . }\textsl{In the Euclidean
space }$L(T)$\textsl{\ one has }

\textsl{i)} $\Lambda =\exp ^{-1}(\mathcal{Z}(G))$\textsl{; ii)} $\Lambda
_{r}\subseteq \Lambda _{e}\subseteq \Lambda $\textsl{, }

\noindent \textsl{where in ii), the first equality\ holds if and only if }$G$
\textsl{is }$1$\textsl{--connected,} \textsl{and} \textsl{the second
equality holds if and only if }$\mathcal{Z}(G)=\{0\}$.$\square $

\bigskip

For a simple Lie group $G$ the quotient group $\Lambda /\Lambda _{e}$ is
always finite (see \cite[p.68]{[Hu]}). As a result we can introduce the 
\textsl{deficiency function }on the weight lattice

\begin{enumerate}
\item[(2.3)] $\kappa :\Lambda \rightarrow \mathbb{Z}$, $x\rightarrow \kappa
_{x}$,
\end{enumerate}

\noindent by letting $\kappa _{x}$ be the least positive integer so that $%
\kappa _{x}x\in \Lambda _{e}$, $x\in \Lambda $. This function provides us
with a partition $\Omega =\Omega _{1}\sqcup \Omega _{2}$ with

\begin{quote}
$\Omega _{1}=\{\omega \in \Omega \mid \kappa _{\omega }=1\}$, $\Omega
_{2}=\{\omega \in \Omega \mid \kappa _{\omega }\geq 2\}$.
\end{quote}

\noindent \textbf{Example 2.6.} Let $G$ be a simple Lie group.

If $\mathcal{Z}(G)=\{0\}$ we get from $\Lambda _{e}=\Lambda $ by Lemma 2.5
that $\Omega =\Omega _{1}$.

If $G$ is $1$--connected with $\mathcal{Z}(G)\neq \{0\}$, we have $\Lambda
_{e}=\Lambda _{r}$ by Lemma 2.5. From the expressions of the fundamental
dominant weights\textsl{\ }by simple roots in \cite[p.69]{[Hu]} one
determines the subset $\Omega _{1}$, consequently $\Omega _{2}$, as that
tabulated below

\begin{center}
\begin{tabular}{l|llllll}
\hline\hline
$G$ & $A_{n}$ & $Sp(n)$ & $Spin(2n+1)$ & $Spin(2n)$ & $E_{6}$ & $E_{7}$ \\ 
\hline
$\Omega _{1}$ & $\emptyset $ & $\{\omega _{i}\}_{i<n}$ & $\{\omega
_{i}\}_{i=2k}$ & $\{\omega _{i}\}_{i=2k\text{ }\leq n-2,}$ & $\{\omega
_{i}\}_{i=2,4}$ & $\{\omega _{i}\}_{i=1,3,4,6}$ \\ \hline\hline
\end{tabular}%
.$\square $
\end{center}

The set $\Delta (x_{0})$ of simple roots is a basis for both $L(T)$ and the
root lattice $\Lambda _{r}$. Using this basis a partial order $\prec $ on $%
L(T)$ (hence on $\Phi \subset L(T)$) can be introduced by the following rule:

\begin{center}
$v\prec u$\textsl{\ if and only if the difference }$u-v$\textsl{\ is a sum
of elements of }$\Delta (x_{0})$\textsl{.}
\end{center}

\noindent As in \cite[p.67]{[Hu]} we put $\Lambda ^{+}:=\Lambda \cap 
\mathcal{F}(x_{0})$. An element $\omega \in \Lambda ^{+}$ is called \textsl{%
minimal }if $\omega \succ \omega ^{\prime }\in \Lambda ^{+}$ implies that $%
\omega =\omega ^{\prime }$.

\bigskip

\noindent \textbf{Lemma 2.7.} \textsl{Let }$G$ \textsl{be an }$1$\textsl{%
--connected simple Lie group, and\ let }$\Pi _{G}\subset \Lambda ^{+}$ 
\textsl{be the subset of all non--zero minimal weights. Then\ }$\Pi
_{G}\subset \Omega _{2}$\textsl{. Moreover,}

\textsl{i) for} \textsl{each} $\omega \in \Omega _{2}$ \textsl{there is
precisely one} \textsl{weight} $\omega ^{\prime }\in \Pi _{G}$ \textsl{so
that }$\omega \succ \omega ^{\prime }$\textsl{;}

\textsl{ii) the set of all non--trivial elements in }$\mathcal{Z}(G)$\textsl{%
\ are given without repetition by }$\{\exp (\omega )\in \mathcal{Z}(G)\mid
\omega \in \Pi _{G}\}$\textsl{.}

\noindent \textbf{Proof.} See \cite[P.92]{[Hu]}. $\square $

\bigskip

In view of Lemma 2.7 we can introduce a \textsl{retraction} $r:\Omega
_{2}\rightarrow \Pi _{G}$ and an \textsl{involution} $\tau :\Pi
_{G}\rightarrow \Pi _{G}$ respectively by the rules:

a) $\omega \succ r(\omega )\in \Pi _{G}$, $\omega \in \Omega _{2}$ (by i) of
Lemma 2.7);

b) $\tau (\omega )+\omega \in \Lambda _{e}$, $\omega \in \Pi _{G}$ (by ii)
of Lemma 2.7).

\noindent Alternatively, the element $\tau (\omega )$ is characterized by
the relation

\begin{quote}
$\exp (\omega )\exp (\tau (\omega ))=e$.
\end{quote}

\noindent \textbf{Example 2.8. }Assume that\textbf{\ }$G$ is simple and $1$%
--connected with $\mathcal{Z}(G)\neq \{0\}$.

a) The set $\Pi _{G}$ of minimal weights is given by (see \cite[P.92]{[Hu]}):

\begin{center}
\begin{tabular}{l|l|l|l|l|l|l}
\hline\hline
$G$ & $SU(n)$ & $Sp(n)$ & $Spin(2n+1)$ & $Spin(2n)$ & $E_{6}$ & $E_{7}$ \\ 
\hline
$\Pi _{G}$ & $\left\{ \omega _{i}\right\} _{1\leq i\leq n}$ & $\left\{
\omega _{n}\right\} $ & $\left\{ \omega _{1}\right\} $ & $\left\{ \omega
_{1},\omega _{n-1},\omega _{n}\right\} $ & $\left\{ \omega _{1},\omega
_{6}\right\} $ & $\left\{ \omega _{7}\right\} $ \\ \hline\hline
\end{tabular}%
;
\end{center}

b) The set $\Omega _{2}$, as well as the composition $\tau \circ r:\Omega
_{2}\rightarrow \Pi _{G}$, is given by

\begin{center}
\begin{tabular}{l|l|l|l|l|l}
\hline\hline
$G$ & $SU(n)$ & $Sp(n)$ & $Spin(2n+1)$ & $E_{6}$ & $E_{7}$ \\ \hline
$\Omega _{2}$ & ${\small \{\omega }_{{\small k}}{\small \}}_{{\small 1\leq
k\leq n}}$ & ${\small \{\omega }_{{\small n}}{\small \}}$ & ${\small %
\{\omega }_{{\small 2k+1}}{\small \}}$ & ${\small \{\omega }_{{\small k}}%
{\small \}}_{{\small k=1,3,5,6}}$ & ${\small \{\omega }_{k}{\small \}}_{%
{\small k=2,5,7}}$ \\ \hline
${\small \tau \circ r(\omega }_{{\small k}}{\small )}$ & $\omega _{n+1-k}$ & 
$\omega _{n}$ & $\omega _{1}$ & 
\begin{tabular}{l}
${\small \omega }_{{\small 6}}\text{ }${\small for }${\small k=1,5}$ \\ 
${\small \omega }_{{\small 1}}$ {\small for }${\small k=3,6}$%
\end{tabular}
& $\omega _{7}$ \\ \hline\hline
\end{tabular}
\end{center}

\noindent and for $G=Spin(2n)$, by $\Omega _{2}=\{\omega _{k},\omega
_{n-1},\omega _{n}\mid k\leq n-2$ odd$\}$,

\begin{center}
$\tau \circ r(\omega _{k})=\left\{ 
\begin{tabular}{l}
$\omega _{1}\text{ if }k\leq n-2\text{;}$ \\ 
$\omega _{n}$ if either $n$ is odd, $k=n-1$, or $n$ is even, $k=n$; \\ 
$\omega _{n-1}$ if either $n$ is odd, $k=n$, or $n$ is even, $k=n-1$.$%
\square $%
\end{tabular}%
\right. $
\end{center}

\section{Computation in the fundamental Weyl cell}

For a simple Lie group $G$ elements in the root system\textsl{\ }$\Phi $ has
at most two lengths. Let $\beta $\textsl{\ }be the \textsl{maximal short root%
} relative to the partial order $\prec $ on the set $\Phi ^{+}$ of positive
roots \cite[p.55]{[Hu]}. The\textsl{\ fundamental Weyl cell} is the simplex
in the Weyl chamber $\mathcal{F}(x_{0})$ defined by

\begin{quote}
$\Delta =\{u\in \mathcal{F}(x_{0})\mid \beta ^{\ast }(u)\leq 1\}$.
\end{quote}

\noindent Let $d:T\times T\rightarrow \mathbb{R}$\ be the distance function
on $T$. It is well known that

\bigskip

\noindent \textbf{Lemma 3.1 (\cite{[C], [Cr]}).} \textsl{Let }$G$\textsl{\
be a simple Lie group. Then}

\textsl{i) the equation }$d(e,\exp (u))=\left\Vert u\right\Vert $\textsl{\
holds if and only if }$\left\Vert u\right\Vert \leq \left\Vert
u-v\right\Vert $ \textsl{holds} \textsl{for all }$v\in \Lambda _{e}$\textsl{;%
}

\textsl{ii)} \textsl{if }$G$\textsl{\ is }$1$\textsl{--connected, then }$%
u\in \Delta $ \textsl{implies that} $d(e,\exp (u))=\left\Vert u\right\Vert $%
\textsl{.}$\square $

\bigskip

It is well known that every element $x\in G$ is conjugate under $G$ to an
element of the form $\exp (u)\in G$ with $u\in \Delta $ and $d(e,\exp
(u))=\left\Vert u\right\Vert $. Moreover, if $x\in Fix(\gamma )$ then $2u\in
\Lambda _{e}$. This implies that

\bigskip

\noindent \textbf{Lemma 3.2.} \textsl{For a simple Lie group }$G$\textsl{\
with fundamental Weyl cell }$\Delta $\textsl{\ set}

\begin{enumerate}
\item[(3.1)] $\mathcal{K}_{G}=\{u\in \Delta \mid 2u\in \Lambda _{e}$,\ $%
d(e,\exp (u))=\left\Vert u\right\Vert \}$.
\end{enumerate}

\noindent \textsl{Then}

\begin{enumerate}
\item[(3.2)] $Fix(\gamma )=\{e\}\bigcup\limits_{u\in \mathcal{K}_{G}}M_{\exp
(u)}$ \textsl{with }$d(e,x)=\left\Vert u\right\Vert $\textsl{\ for all }$%
x\in M_{\exp (u)}$\textsl{.}$\square $
\end{enumerate}

Comparing (3.2) with (1.2) we emphasis that the decomposition (3.2) on $%
Fix(\gamma )$ may not be disjoint, as overlap like $M_{\exp (u)}=M_{\exp
(v)} $ may occur for some $u,v\in \mathcal{K}_{G}$ with $u\neq v$. However,
based on the relation (3.2) our approach to $Fix(\gamma )$ consists of three
steps:

i) find a general expression for elements in $\mathcal{K}_{G}$;

ii) specify a subset $\mathcal{F}_{G}\subseteq \mathcal{K}_{G}$ so that the
relation (3.2) can be refined as $Fix(\gamma )=\{e\}\coprod\limits_{u\in 
\mathcal{F}_{G}}M_{\exp (u)}$;

iii) decide the isomorphism types of $M_{\exp (u)}$ for all $u\in \mathcal{F}%
_{G}$.

\noindent In this section we accomplish step i) in the next result.

\bigskip

\noindent \textbf{Lemma 3.3. }\textsl{Let} $G$ \textsl{be a simple Lie group.%
} \textsl{Then }$u\in \mathcal{K}_{G}$ \textsl{implies that}

\begin{enumerate}
\item[(3.3)] $u=\left\{ 
\begin{tabular}{l}
$\frac{1}{2}\omega _{k}$ \textsl{for some }$\omega _{k}\in \Omega _{1}$, \\ 
$\frac{1}{2}(\omega _{k}+\tau \circ r(\omega _{k}))$ \textsl{for some} $%
\omega _{k}\in \Omega _{2}$\textsl{.}%
\end{tabular}%
\right. $
\end{enumerate}

\noindent \textbf{Proof.} For an $u\in \mathcal{K}_{G}$ we get from $u\in
\Delta $ and $2u\in \Lambda _{e}$ that

\begin{quote}
$u=\lambda _{k_{1}}\omega _{k_{1}}+\cdots +\lambda _{k_{t}}\omega _{k_{t}}$
with $\lambda _{k_{s}}>0$ and $2\lambda _{k_{s}}\in \mathbb{Z}$
\end{quote}

\noindent by Lemma 2.4, where $\{k_{1},\cdots ,k_{t}\}\subseteq \{1,\cdots
,n\}$. This implies that

\begin{enumerate}
\item[(3.4)] $2u-\omega _{k_{1}}=a\in \Lambda ^{+}$.
\end{enumerate}

\noindent The formula (3.3) will be deduced from the second constraint $%
d(e,\exp (u))=\left\Vert u\right\Vert $ on $u\in \mathcal{K}_{G}$ in (3.1).

If $\omega _{k_{1}}\in \Omega _{1}$ then $\omega _{k_{1}}\in \Lambda _{e}$
implies that $\left\Vert u\right\Vert \leq \left\Vert u-\omega
_{k_{1}}\right\Vert $ by i) of Lemma 3.1. That is

\begin{quote}
$\left\Vert \frac{1}{2}a+\frac{1}{2}\omega _{k_{1}}\right\Vert ^{2}\leq
\left\Vert \frac{1}{2}a-\frac{1}{2}\omega _{k_{1}}\right\Vert ^{2}$
\end{quote}

\noindent by (3.4). However, since $(\omega _{i},\omega _{j})>0$ by iii) of
Lemma 2.4 and since $a\in \Lambda ^{+}$, this is possible if and only if $%
a=0 $. That is

\begin{enumerate}
\item[(3.5)] $u=\frac{1}{2}\omega _{k_{1}}$.
\end{enumerate}

If $\omega _{k_{1}}\in \Omega _{2}$ we have $a\in \Lambda ^{+}$ but $a\notin
\Lambda _{e}$ by (3.4). According to Lemma 2.7 there is precisely one weight 
$\omega _{s}\in \Pi _{G}$ so that $2u-\omega _{k_{1}}=a=\omega _{s}+b$ with $%
b$ a sum of elements of $\Delta (x_{0})$. From $2u,b\in \Lambda _{e}$ we
find that $\omega _{k_{1}}+\omega _{s}\in \Lambda _{e}$ and therefore $%
\left\Vert u\right\Vert \leq \left\Vert u-\omega _{k_{1}}-\omega
_{s}\right\Vert $ by i) of Lemma 3.1. That is

\begin{quote}
$\left\Vert \frac{1}{2}(\omega _{k_{1}}+\omega _{s})+\frac{1}{2}b\right\Vert
^{2}\leq \left\Vert \frac{1}{2}(\omega _{k_{1}}+\omega _{s})-\frac{1}{2}%
b\right\Vert ^{2}$.
\end{quote}

\noindent Again, since $(\omega _{i},\omega _{j})>0$ by iii) of Lemma 2.4
and since $b\in \Lambda ^{+}$, this is possible if and only if $b=0$. We
obtain from (3.4) that $u=\frac{1}{2}(\omega _{k_{1}}+\omega _{s})$, $\omega
_{s}\in \Pi _{G}$. Furthermore, from the calculation

\begin{quote}
$e=\exp (2u)=\exp (\omega _{k_{1}})\exp (\omega _{s})=\exp (r(\omega
_{k_{1}}))\exp (\omega _{s})$
\end{quote}

\noindent (since $\omega _{k_{1}}\succ r(\omega _{k_{1}})\in \Pi _{G}$) as
well as the definition of $\tau $ we get $\omega _{s}=$ $\tau \circ r(\omega
_{k_{1}})$. This shows that

\begin{enumerate}
\item[(3.6)] $u=\frac{1}{2}(\omega _{k_{1}}+\tau \circ r(\omega _{k_{1}}))$
with $\omega _{k_{1}}\in \Omega _{2}$.
\end{enumerate}

\noindent The proof of (3.3) has now been completed by (3.5) and (3.6).$%
\square $

\section{Proof of Theorem 1.1}

Assume that $G$ is an exceptional Lie group and the expression of its
maximal short root $\beta $ in term of the simple roots is\ $\beta
=m_{1}\alpha _{1}+\cdots +m_{n}\alpha _{n}$ (see \cite[p.66]{[Hu]}). By the
definition (2.2) of the fundamental dominant weights

\begin{enumerate}
\item[(4.1)] $\beta ^{\ast }(\frac{\omega _{i}}{2})=\frac{m_{i}\left\Vert
\alpha _{i}\right\Vert ^{2}}{2\left\Vert \beta \right\Vert ^{2}}$; $\quad
\beta ^{\ast }(\frac{\omega _{i}+\omega _{j}}{2})=\frac{m_{i}\left\Vert
\alpha _{i}\right\Vert ^{2}}{2\left\Vert \beta \right\Vert ^{2}}+\frac{%
m_{j}\left\Vert \alpha _{j}\right\Vert ^{2}}{2\left\Vert \beta \right\Vert
^{2}}$.
\end{enumerate}

\noindent Let $\mathcal{K}_{G}^{\prime }\subset \Delta $ be the subset of
the vectors $u$ satisfying (3.3). Combining (3.3) and (4.1), together with
computations in Examples 2.6 and 2.8, one determines the set $\mathcal{K}%
_{G}^{\prime }$ for each exceptional $G$, as that presented in the second
column of Tables 2 below.

In \cite{[DL]} an explicit procedure to calculate the isomorphism type of
the centralizer $C_{\exp (u)}\subset G$ in term of $u\in \Delta $ is
obtained. As applications those centralizers $C_{\exp (u)}$ with $u\in 
\mathcal{K}_{G}^{\prime }$ are determined and presented in the third column
of Table 2 (see also in \cite[Theorem 4.4, Theorem 4.6]{[DL]}).

In general $\mathcal{K}_{G}\subseteq \mathcal{K}_{G}^{\prime }$ by Lemma
3.3. However, the centralizers $C_{\exp (u)}$ recorded in Table 2 are useful
for us to specify the desired subset $\mathcal{K}_{G}$ from $\mathcal{K}%
_{G}^{\prime }$. To explain this we observe that if $u\in \mathcal{K}%
_{G}^{\prime }$ is a vector with $u\notin \mathcal{K}_{G}$, then $d(e,\exp
(u))<\left\Vert u\right\Vert $ implies that there exists a vector $v\in L(T)$
satisfying

\begin{quote}
$\exp (v)=\exp (u)$ and $d(e,\exp (v))=\left\Vert v\right\Vert $.
\end{quote}

\noindent Take a Weyl group element $w\in W$ so that $v^{\prime }=w(v)\in 
\mathcal{F}(x_{0})$. The relations $2v^{\prime }\in \Lambda _{e}$ and $%
d(e,\exp (v^{\prime }))=\left\Vert v^{\prime }\right\Vert $ indicate that $%
v^{\prime }\in \mathcal{K}_{G}$ by (3.1). In particular we have

\bigskip

\noindent \textbf{Lemma 4.1.} If $u\in \mathcal{K}_{G}^{\prime }\backslash 
\mathcal{K}_{G}$ then there exists an element $v\prime \in \mathcal{K}_{G}$
so that

\begin{quote}
$C_{\exp (u)}\cong C_{\exp (v^{\prime })}$ and $\left\Vert v^{\prime
}\right\Vert <\left\Vert u\right\Vert $.$\square $
\end{quote}

\begin{center}
\begin{tabular}{l|l|l}
\hline\hline
$G$ & $\mathcal{K}_{G}^{\prime }$ & the centralizer $C_{\exp (u)}$ with $%
u\in \mathcal{K}_{G}^{\prime }$ \\ \hline
$G_{2}$ & $\{\frac{1}{2}\omega _{{\small 1}}\}$ & ${\small SO(4)}$ \\ 
$F_{4}$ & $\{\frac{1}{2}\omega _{{\small k}}\}_{{\small k=1,4}}$ & ${\small %
Spin(9),}\frac{{\small Sp(3)\times Sp(1)}}{{\small Z}_{2}}$ \\ 
$E_{6}$ & $\{\frac{{\small \omega }_{{\small 2}}}{2},\frac{\omega
_{1}+\omega _{6}}{2}\}$ & $\frac{{\small SU(2)\times SU(6)}}{{\small Z}_{2}}$%
{\small ,} $\frac{{\small Spin(10)\times S}^{{\small 1}}}{{\small Z}_{4}}$
\\ 
$E_{6}^{\ast }$ & ${\small \{}\frac{1}{2}{\small \omega }_{{\small k}}%
{\small \}}_{{\small k=1,2,3,5},{\small 6}}$ & $\frac{{\small Spin(10)\times
S}^{{\small 1}}}{{\small Z}_{4}}$ {\small for }${\small k=1,6}${\small ; }$%
\frac{{\small SU(2)\times }\frac{{\small SU(6)}}{Z_{3}}}{{\small Z}_{2}}$ 
{\small for }${\small k=2,3,5}$ \\ 
$E_{7}$ & $\{\frac{{\small \omega }_{{\small k}}}{2},\omega _{{\small 7}}\}_{%
{\small k=1,6}}$ & $\frac{{\small Spin(12)\times SU(2)}}{{\small Z}_{2}}$%
{\small , }$\frac{{\small Spin(12)\times SU(2)}}{{\small Z}_{2}}${\small , }$%
{\small E}_{{\small 7}}$ \\ 
$E_{7}^{\ast }$ & ${\small \{}\frac{1}{2}{\small \omega }_{{\small k}}%
{\small \}}_{{\small k=1,2,6,7}}$ & $\frac{Ss(12){\small \times SU(2)}}{%
{\small Z}_{2}}${\small ,} $[\frac{{\small SU(8)}}{{\small Z}_{4}}]^{2}$%
{\small , }$\frac{Ss(12){\small \times SU(2)}}{{\small Z}_{2}}$, $[\frac{%
{\small E}_{6}{\small \times S}^{1}}{{\small Z}_{3}}]^{2}$ \\ 
$E_{8}$ & $\{\frac{1}{2}\omega _{{\small k}}\}_{{\small k=1,8}}$ & ${\small %
Ss(16)}${\small , }$\frac{{\small E}_{{\small 7}}{\small \times SU(2)}}{%
{\small Z}_{2}}$ \\ \hline\hline
\end{tabular}

{\small Table 2. the set }$\mathcal{K}_{G}^{\prime }${\small \ as well as
the centralizers }${\small C}_{\exp (u)}${\small \ with }$u\in \mathcal{K}%
_{G}^{\prime }$.
\end{center}

\noindent \textbf{Proof of Theorem 1.1.} The proof will be divided into two
cases, depending on whether $G$ is $1$--connected. Concerning the use of
Lemma 4.1 in the forthcoming arguments, we note that in view of the
presentation of the fundamental dominant weights with respect to appropriate
Euclidean coordinates $\{\varepsilon _{1},\cdots ,\varepsilon _{m}\}$ on $%
L(T)$ in the standard reference \cite[p.265-277]{[B]}, the length $%
\left\Vert u\right\Vert $ for a vector $u\in \mathcal{K}_{G}^{\prime }$ can
be easily evaluated.

\textbf{Case I.} $G=G_{2},F_{4},E_{6},E_{7},E_{8}$. Since $G$ is $1$%
--connected, we have

\begin{enumerate}
\item[(4.2)] $\mathcal{K}_{G}=\mathcal{K}_{G}^{\prime }$
\end{enumerate}

\noindent by ii) of Lemma 3.1. Consequently, it follows from (3.2) that

\begin{enumerate}
\item[(4.3)] $Fix(\gamma )=\{e\}\bigcup\limits_{u\in \mathcal{K}_{G}}M_{\exp
(u)}$ with $M_{\exp (u)}\cong G/C_{\exp (u)}$.
\end{enumerate}

Assume in the decomposition (4.3) on $Fix(\gamma )$ that the relation $%
M_{\exp (u)}=M_{\exp (v)}$ holds for some $u,v\in \mathcal{K}_{G}$ with $%
u\neq v$. By (3.2) one has

\begin{quote}
i) $C_{\exp (u)}\cong C_{\exp (v)}$ and ii) $\left\Vert u\right\Vert
=\left\Vert v\right\Vert $.
\end{quote}

\noindent However, in view of the groups $C_{\exp (u)}$ presented in Table 2
the only possibility for i) to hold is when $G=E_{7}$ and $(u,v)=(\frac{%
{\small \omega }_{{\small 1}}}{2},\frac{{\small \omega }_{{\small 6}}}{2})$,
but in this case the calculation

\begin{quote}
$\left\Vert \frac{{\small \omega }_{{\small 1}}}{2}\right\Vert =\frac{1}{%
\sqrt{2}}<\left\Vert \frac{{\small \omega }_{{\small 6}}}{2}\right\Vert =1$
(see \cite[p.280]{[B]}).
\end{quote}

\noindent shows that the relation ii) does not hold. Summarizing, taking $%
\mathcal{F}_{G}=\mathcal{K}_{G}$ then the decomposition (1.2) on $Fix(\gamma
)$ is given by (4.3), and the proof of Theorem 1.1 is completed by the
corresponding items in Table 2.

\textbf{Case II.} $G=E_{6}^{\ast },E_{7}^{\ast }$. This case is slightly
delicate because, instead of the equality (4.2) one has $\mathcal{K}%
_{G}\subseteq \mathcal{K}_{G}^{\prime }$ by Lemma 3.3. Nevertheless, granted
with Lemma 4.1 and results in Table 2 we shall show that

\begin{enumerate}
\item[(4.4)] $\mathcal{K}_{G}=\left\{ 
\begin{tabular}{l}
${\small \{}\frac{{\small \omega }_{{\small k}}}{2}{\small \}}_{{\small %
k=1,2,6}}$ if $G=E_{6}^{\ast }$ \\ 
${\small \{}\frac{{\small \omega }_{{\small k}}}{2}{\small \}}_{{\small %
k=1,2,7}}$ if $G=E_{7}^{\ast }$%
\end{tabular}%
\right. ;$

\item[(4.5)] for $u,v\in \mathcal{K}_{G}$ the overlap $M_{\exp (u)}=M_{\exp
(v)}$ (see (3.2)) happens if and only if $G=E_{6}^{\ast }$ and $(u,v)=(\frac{%
{\small \omega }_{1}}{2},\frac{{\small \omega }_{{\small 6}}}{2})$.
\end{enumerate}

\noindent Consequently, setting

\begin{quote}
$\mathcal{F}_{G}=\left\{ 
\begin{tabular}{l}
${\small \{}\frac{{\small \omega }_{{\small k}}}{2}{\small \}}_{{\small k=1,2%
}}$ if $G=E_{6}^{\ast }$ \\ 
${\small \{}\frac{{\small \omega }_{{\small k}}}{2}{\small \}}_{{\small %
k=1,2,7}}$ if $G=E_{7}^{\ast }$%
\end{tabular}%
\right. $
\end{quote}

\noindent the proof of Theorem 1.1 for this case is completed by (4.4) and
(4.5), and the relevant items in Table 2.

For $G=E_{7}^{\ast }$ we have $\mathcal{K}_{E_{7}^{\ast }}^{\prime }={\small %
\{}\frac{1}{2}{\small \omega }_{{\small k}}{\small \}}_{{\small k=1,2,6,7}}$
by Table 2. Since $\omega _{7}\in \Lambda _{e}$ by Example 2.6 and since

\begin{quote}
$\frac{1}{\sqrt{2}}=\left\Vert \omega _{7}-\frac{1}{2}\omega _{6}\right\Vert
<\left\Vert \frac{1}{2}\omega _{6}\right\Vert =1$ (see \cite[p.280]{[B]})
\end{quote}

\noindent we have $\frac{1}{2}\omega _{6}\notin \mathcal{K}_{E_{7}^{\ast }}$
by i) of Lemma 3.1. The proof of (4.4) for $G=E_{7}^{\ast }$ is done by
Lemma 4.1, together with the groups $C_{\exp (u)}$ with $u\in {\small \{}%
\frac{1}{2}{\small \omega }_{{\small k}}{\small \}}_{{\small k=1,2,7}}$
given in Table 2.

Similarly, for $G=E_{6}^{\ast }$ we have $\omega _{1},\omega _{6}\in \Lambda
_{e}$ by Example 2.6, but

\begin{quote}
$\frac{1}{\sqrt{2}}=\left\Vert \omega _{1}-\frac{1}{2}\omega _{3}\right\Vert
<\left\Vert \frac{1}{2}\omega _{3}\right\Vert =\sqrt{\frac{5}{6}}$;

$\frac{1}{\sqrt{2}}=\left\Vert \omega _{6}-\frac{1}{2}\omega _{5}\right\Vert
<\left\Vert \frac{1}{2}\omega _{5}\right\Vert =\sqrt{\frac{5}{6}}$, see \cite%
[p.276]{[B]}.
\end{quote}

\noindent We get $\frac{1}{2}\omega _{3},\frac{1}{2}\omega _{5}\notin 
\mathcal{K}_{E_{6}^{\ast }}$ from i) of Lemma 3.3. The proof of (4.4) for $%
G=E_{6}^{\ast }$ is done by Lemma 4.1, together with the groups $C_{\exp
(u)} $ with $u\in {\small \{}\frac{1}{2}{\small \omega }_{{\small k}}{\small %
\}}_{{\small k=1,2,6}}$ given in Table 2.

For (4.5) assume that in the decomposition (3.2) on $Fix(\gamma )$ the
relation $M_{\exp (u)}=M_{\exp (v)}$ holds for some $u,v\in \mathcal{K}_{G}$
with $u\neq v$. By (3.2) one has

\begin{quote}
i) $C_{\exp (u)}\cong C_{\exp (v)}$; ii) $\left\Vert u\right\Vert
=\left\Vert v\right\Vert $.
\end{quote}

\noindent In view of the groups $C_{\exp (u)}$ with $u\in \mathcal{K}_{G}$
presented in the last column of Table 2 the only possibility for both i) and
ii) to hold is when $G=E_{6}^{\ast }$ and $(u,v)=(\frac{{\small \omega }_{%
{\small 1}}}{2},\frac{{\small \omega }_{{\small 6}}}{2})$. Let $w_{0}$ be
the unique longest element of the Weyl group of $E_{6}$ \cite[p.171]{[B]}.
Then

\begin{quote}
$w_{0}(\frac{1}{2}\omega _{{\small 1}})=-\frac{1}{2}\omega _{{\small 6}}$
(see \cite[p.276]{[B]}).
\end{quote}

\noindent This implies that $M_{\exp (\frac{{\small \omega }_{{\small 1}}}{2}%
)}=M_{\exp (-\frac{1}{2}\omega _{{\small 6}})}$. We get (4.5) from the
general relation $\exp (-\frac{1}{2}u)=\exp (\frac{1}{2}u)$, $u\in \Lambda $%
, which holds in all groups $G$ with trivial center. This completes the
proof.$\square $

\bigskip

\noindent \textbf{Remark 4.2. }For a vector $u\in \Delta $ let $C_{\exp
(u)}^{0}$ be the identity component of the centralizer $C_{\exp (u)}$.
Indeed, for $G=E_{6}^{\ast }$ or $E_{7}^{\ast }$ the main result in \cite[%
Theorem 3.7]{[DL]} is applicable to determine the isomorphism type of $%
C_{\exp (u)}^{0}$ instead of the whole group $C_{\exp (u)}$. Therefore,
additional explanation for the groups $C_{\exp (u)}$ corresponding to $%
G=E_{6}^{\ast }$ or $E_{7}^{\ast }$ in Table 2 is requested.

In general, let $p:G^{\symbol{126}}\rightarrow G$ be the universal covering
of a simple Lie group $G$ and let $T^{\symbol{126}}\subset G^{\symbol{126}}$
be the maximal torus of $G^{\symbol{126}}$ corresponding to $T$ in $G$. With
respect to the standard identification $L(G^{\symbol{126}})=L(G)$ (resp. $%
L(T^{\symbol{126}})=L(T)$) the exponential map $\exp $ of $G$ factors
through that $\exp ^{\symbol{126}}$ of $G^{\symbol{126}}$ in the fashion

\begin{quote}
$\exp =p\circ \exp ^{\symbol{126}}:L(G^{\symbol{126}})\rightarrow G^{\symbol{%
126}}\rightarrow G$ (resp. $L(T^{\symbol{126}})\rightarrow T^{\symbol{126}%
}\rightarrow T$).
\end{quote}

\noindent Since for $u\in \mathcal{F}_{G}$ the subspace $M_{\exp ^{\symbol{%
126}}(u)}$ of $G^{\symbol{126}}$ is $1$--connected\textbf{\ }\cite[Corollary
3.4, p.101]{[Bo]}, $p$ restricts to a universal covering $p_{u}:M_{\exp ^{%
\symbol{126}}(u)}\rightarrow M_{\exp (u)}$.

On the other hand, as the set of all non--trivial covering transformations
of $p$ are in one to one correspondence with the set $\Pi _{G^{\symbol{126}%
}} $ of minimal weights (see Example 2.8) in the fashion

\begin{quote}
$\widetilde{g}\rightarrow \exp ^{\symbol{126}}(\omega _{s})\cdot \widetilde{g%
}$, $\widetilde{g}\in G^{\symbol{126}}$, $\omega _{s}\in \Pi _{G^{\symbol{126%
}}}$,
\end{quote}

\noindent the set $\Pi _{u}$ of nontrivial covering transformations of $%
p_{u} $ can be shown to be

\begin{enumerate}
\item[(4.6)] $\Pi _{u}=\{\omega _{s}\in \Pi _{G^{\symbol{126}}}\mid $ $%
\omega _{s}+u-w(u)\in \Lambda _{r}$ for some $w\in W\}$,
\end{enumerate}

\noindent where $\Lambda _{r}$ is the root lattice of $G^{\symbol{126}}$.
Based on this formula a direct calculation in the vector space space $L(T^{%
\symbol{126}})$ shows that

\begin{enumerate}
\item[(4.7)] $\Pi _{u}=\left\{ 
\begin{tabular}{l}
$\{\omega _{7}\}$ for $G=E_{7}^{\ast }$ and $u\in {\small \{}\frac{{\small %
\omega }_{{\small k}}}{2}{\small \}}_{{\small k=2,7}}$ \\ 
$\emptyset $ otherwise.%
\end{tabular}%
\right. $.
\end{enumerate}

\noindent Consequently

\begin{quote}
$C_{\exp (u)}=\left\{ 
\begin{tabular}{l}
$\lbrack C_{\exp (u)}^{0}]^{2}$ for $G=E_{7}^{\ast }$ and $u\in {\small \{}%
\frac{{\small \omega }_{{\small k}}}{2}{\small \}}_{{\small k=2,7}}$ \\ 
$C_{\exp (u)}^{0}$ otherwise.%
\end{tabular}%
\right. $
\end{quote}

\noindent This justify the groups $C_{\exp (u)}$ corresponding to $%
G=E_{6}^{\ast }$ or $E_{7}^{\ast }$ in Table 2.$\square $

\bigskip

\noindent \textbf{Remark 4.3.} With the preliminary data for $%
G=SU(n+1),Sp(n),$ $Spin(n+2)$, $n\geq 2$, recorded in Example 2.8, one can
obtain the fixed set $Fix(\gamma )$ for the simple Lie groups of the
classical types (see \cite{[F]}) by the same argument as that used to
establish Theorem 1.1.$\square $

\bigskip

\noindent \textbf{Remark 4.4.} It is clear that $x\in Fix(\gamma )$\textsl{\ 
}implies that $x^{2}=e$. Consequently, the map $\sigma :G\rightarrow G$ by $%
\sigma (g)=xgx^{-1}$ is an involutive automorphism of $G$ with fixed
subgroup $C_{x}$, the centralizer at $x$. This indicates that the orbit
spaces $M_{\exp (u)}$ in the decomposition (1.2) are all\textsl{\ global
Riemannian symmetric spaces} of $G$ in the sense of E. Cartan. However, the
existing theory of symmetric spaces \cite{[He],[CN1],[N],[Y1],[Y2],[Y3]}
does not constitute a solution to our problem for the following reasons:

i) not every symmetric space of $G$ can appear as a component of $Fix(\gamma
)$;

ii) if a symmetric space of $G$ happens to be a component of $Fix(\gamma )$,
it may occur twice (see in (1.2) for the case $G=E_{7}$);

iii) in view of the relation $M_{\exp (u)}=G/C_{\exp (u)}$ a complete
characterization of the symmetric space $M_{\exp (u)}$ amounts to the
determination of the centralizer $C_{\exp (u)}$, which is a delicate issue
absent in the classical theory of Lie groups \cite{[He]}, and has recently
been made explicit in our paper \cite{[DL]}.

\noindent As a witness of i)--iii), for the exceptional Lie groups Nagano 
\cite{[N]} stated the list of all the symmetric spaces which are connected
components of\noindent\ $Fix(\gamma )$ without specifying their embedding in 
$G$. He did not write a proof in his other papers, although he promised to
do so in \cite{[N]}. In addition, the cases $G=E_{6}^{\ast }$ or $%
E_{7}^{\ast }$ were not considered in the papers \cite{[Y1],[Y2],[Y3]}.

Summarizing, without resorting to the theory of symmetric spaces and by a
unified approach, we have enumerated all the symmetric spaces of an
exceptional $G$ that are components of $Fix(\gamma )$, and presented
concrete realization of these spaces as the adjoint orbits of $G$.$\square $

\section{Generalities}

Result on $Fix(\gamma )$ for the simple Lie groups (i.e. Theorem 1.1 and 
\cite{[F]}) is fundamental in understanding the general structure of the
fixed set $Fix(\gamma )$ for the inverse involution $\gamma $ on an
arbitrary Lie group $G$. Fairly transparent in our context we have the next
result which indicates how Theorem 1.1 could be extended to general settings.

\bigskip

\noindent \textbf{Corollary 5.1.} \textsl{For any semi--simple Lie group }$G$%
\textsl{\ with a maximal torus }$T$\textsl{, there is a finite subset }$%
\mathcal{F}_{G}\subset L(T)$\textsl{\ so that}

\textsl{i) }$\left\Vert u\right\Vert =d(e,\exp (u))$\textsl{\ for all }$u\in 
\mathcal{F}_{G}$\textsl{;}

\textsl{ii) }$Fix(\gamma )=\{e\}\coprod\limits_{u\in \mathcal{F}_{G}}M_{\exp
(u)}$\textsl{.}

\noindent \textsl{In particular, if }$G=G_{1}\times \cdots \times G_{k}$%
\textsl{\ with all the factor groups }$G_{i}$\textsl{\ exceptional, one can
take }$\mathcal{F}_{G}=\mathcal{F}_{G_{1}}\times \cdots \times \mathcal{F}%
_{G_{k}}$\textsl{\ with }$\mathcal{F}_{G_{i}}$ \textsl{being given by the
second column of Table 1. Consequently, for an }$u=(u_{1},\cdots ,u_{k})\in 
\mathcal{F}_{G}$\textsl{\ with }$u_{i}\in \mathcal{F}_{G_{i}}$\textsl{\ one
has}

\begin{quote}
$M_{\exp (u)}=M_{\exp (u_{1})}\times \cdots \times M_{\exp (u_{k})}$\textsl{.%
}$\square $
\end{quote}

A homomorphism $h:G\rightarrow G^{\prime }$ of two semisimple Lie groups $G$
and $G\prime $ clearly satisfies the relation $h(Fix(\gamma ))\subseteq
Fix(\gamma ^{\prime })$. This indicates that Corollary 5.1 can play a role
in the representation theory of Lie groups. More precisely

\bigskip

\noindent \textbf{Corollary 5.2.} \textsl{A group homomorphism }$%
h:G\rightarrow G^{\prime }$\textsl{\ determines uniquely a correspondence }$%
h^{\circ }:\mathcal{F}_{G}\rightarrow \mathcal{F}_{G^{\prime }}\sqcup \{0\}$%
\textsl{\ so that}

\begin{enumerate}
\item[(5.1)] $h(M_{\exp (u)})\subseteq M_{\exp (h^{\circ }(u))}$\textsl{,} $%
u\in \mathcal{F}_{G}$.
\end{enumerate}

\noindent \textsl{Moreover, if }$h:G\rightarrow G^{\prime }$\textsl{\ is the
inclusion of a totally geodesic subgroup, then}

\begin{enumerate}
\item[(5.2)] $\left\Vert h^{\circ }(u)\right\Vert =\left\Vert u\right\Vert $%
\textsl{\ for all }$u\in \mathcal{F}_{G}$\textsl{.}$\square $
\end{enumerate}

Specifying a subset $\mathcal{F}_{G}\subset L(T)$ with properties (1.2)
amounts to an explicit characterization of the embedding $Fix(\gamma
)\subset G$. Apart from the general fact that the choice of $\mathcal{F}_{G}$
may not be unique, our proof of Theorem 1.1 implies that, if $G$ is $1$%
--connected, there exists a unique set $\mathcal{F}_{G}$ satisfying the
relation $\mathcal{F}_{G}\subset \Delta $. Geometrically

\bigskip

\noindent \textbf{Corollary 5.3. }\textsl{If }$G$\textsl{\ is simple and }$1$%
\textsl{--connected,}\textbf{\ }\textsl{each adjoint orbit }$M_{\exp (u)}$ 
\textsl{in }$Fix(\gamma )$\textsl{\ meets the subspace }$\exp (\Delta )$%
\textsl{\ of }$G$\textsl{\ exactly at one point.}$\square $

\bigskip

\textbf{Acknowledgement.} The authors are grateful to Angela Pasquale for
valuable communications, and in particular, for informing us the works \cite%
{[Y1],[Y2],[Y3]} by I. Yokota.

\end{document}